\documentclass[final,11pt]{article}
\usepackage{latexsym}
\usepackage{multicol}
\usepackage{fancybox}
     \usepackage{xcolor}
\usepackage{amsmath}
\usepackage{cite}

\newtheorem{theorem}{Theorem}

\newtheorem{remark}{Remark}
\newtheorem{example}{Example}[section]
\renewcommand{\baselinestretch}{1.1}
\def\eqsp{\noalign{\vskip 5pt}}
\def\dsp{\displaystyle}

\usepackage{bm}
\usepackage{graphicx}
\usepackage{framed}

\usepackage[all]{xy}

\newcommand{\eqn}[1]{(\ref{#1})}

\newcommand{\eq}{\begin{equation}}
\newcommand{\en}{\end{equation}}
\newcommand{\eqm}{\begin{eqnarray}}
\newcommand{\enm}{\end{eqnarray}}
\newcommand{\eqmno}{\begin{eqnarray*}}
\newcommand{\enmno}{\end{eqnarray*}}
\newcommand{\eql}[1]{\begin{equation}\label{#1}}

\newcommand{\fr}[2]{{\frac{#1}{#2}} }

\newcommand{\goto}{\rightarrow}

\newcommand{\ignore}[1]{}
\def\bc{\begin{center}}
\def\ec{\end{center}}
\def\bi{\begin{itemize}}
\def\ei{\end{itemize}}
\def\be{\begin{enumerate}}
\def\ee{\end{enumerate}}

\def\reals{{{\rm l} \kern -.15em {\rm R} }}

\def\qquad{\quad\quad}

\def\eqsp{\noalign{\vskip 5pt}}
\def\qquad{\quad\quad}

\def\dsp{\displaystyle}
\newcommand{\eqml}[1]{\eql{#1}\begin{array}{rcl}}
\newcommand{\enml}{\end{array}\end{equation}}

\usepackage{multicol}
\usepackage{fancybox}
     \usepackage{xcolor}
  \usepackage{showlabels}
  \usepackage{hyperref}
\usepackage{cite}
 \usepackage{showkeys}

\usepackage{cases}
\usepackage{multirow}
\usepackage{multicol}
\usepackage{arydshln}
\usepackage{amsmath}
\usepackage{amssymb}

\renewcommand{\baselinestretch}{1.1}

\usepackage{bm}
\usepackage{graphicx}
\usepackage{framed}

\usepackage[all]{xy}

\begin{document}

\title{A simple third order compact  finite element method for 1D BVP}

\author{Baiying Dong\thanks{School of Mathematics and Statistics, NingXia University \& School of Mathematics and
Computer Science, NingXia Normal University, China  }
\and Zhilin Li \thanks{CRSC \&
Department of Mathematics, North Carolina State University, Raleigh,
NC 27695-8205, USA}
}

\date{}
\maketitle


\begin{abstract}
A simple third order compact  finite element method  is proposed for one-dimensional Sturm-Liouville boundary value problems.  The key idea is  based on the interpolation error estimate, which can be related to the source term.
Thus, a simple posterior error analysis or a modified basis functions based on original piecewise linear basis function will lead to a third order accurate solution in the $L^2$ norm, and second order in the $H^1$ or the energy norm. Numerical examples have confirmed our analysis.
 \end{abstract}

{\bf keywords:}
compact scheme, finite element method, high order FEM, posterior analysis.



\section{Introduction}

In this paper, we discuss a new third order finite element method for the following Sturm-Liouville boundary value problem
\eqml{1deq}
 \dsp - \nabla \cdot \Big(\beta(x)\nabla u(x) \Big ) + q(x) u(x)=f(x),  \qquad x_l < x < x_r,
\enml
with  linear boundary conditions at $x=a$ and $x=b$, such as Dirichlet, Neumann, and Robin condition.
 To guarantee the well-posedness of the problem, we assume that
 \eqm
   \beta(x)\in L^{\infty}(x_l,x_r) , \quad  \beta(x) \ge \beta_0 >0, \quad q(x)  \in L^{\infty}(x_l,x_r), \quad q(x)\ge 0, \quad f(x)\in L^2(x_l,x_r),
 \enm
 where $\beta_0$ is a constant. With one of Dirichlet boundary condition, that is, the solution is specified at $x=x_r$ or $x=x_r$, then the problem has a unique solution $u(x)\in H^2(x_l,x_r)$ from the Lax-Milgram lemma, see for example, \cite{braess,li:book-FDFEM}.

It is often difficult if not possible to obtain the analytic solution to the problem. In many situations, it is easier to get a numerical solution using modern computers.  Basic finite difference (FD) or finite element (FE) methods have been well developed and understood. In a FD or FE method, a mesh is constructed with a mesh size parameter $h$ so that an approximation to the true solution can be approximated. Commonly used ones are second order accurate meaning that the errors of the computed solution in approximating  the true   is of $O(h^2)$ in some norms, see for example, \cite{Brenner-Scott,morton-mayers}.

Naturally, we hope to have better than second order methods, like third order ($O(h^3)$) or higher without too much additional cost and effort.
There are many high order compact finite difference methods in the literature but  almost none high order compact   finite element  methods.
In this paper, we propose a third order accurate finite element method based on a novel idea which is dependent on existing $P_1$ (using piecewise linear basis functions) finite element method. The new idea is based on the interpolation theory to add a local compact support to the basis function so that the high order term in the error is of $O(h^3)$.



\section{The classical FEM using $P_1$ finite element}

In the classical  finite element method using a piecewise linear function to approximate the solution to \eqn{1deq}, a mesh is generated
\eqm
  x_l = x_0 < x_1 < x_2 < \cdots < x_i < \cdots < x_{n-1} < x_n = x_r,
\enm
where $n$ is a parameter and the mesh parameter is defined as $h = \dsp \max_{1\leq i \leq n} | x_{i+1} - x_i|$. Note that $h\sim O(1/n)$.

In the finite element method, the weak form is below
\eqml{weak}
  & \dsp a(u,v) = L(v), \qquad \forall v(x) \in H^1_0(x_l,x_r), \\ \eqsp
  & \dsp a(u,v) = \int_{x_l}^{x_r} \Big ( \beta(x)u'(x) v'(x) +q(x)u(x)v(x) \Big) dx + R_1(u,v),  \\ \eqsp
  & \dsp  L(v) =  \int_{x_l}^{x_r}   f v dx + R_2(v),
\enml
where $R_1$ and $R_2$  are  contributions from boundary conditions, $a(u,v)$ is a bilinear form, and $L(v)$ is a linear form.
For the homogenous Dirichlet boundary conditions $u(x_l)=u(x_r)=0$,  $R_1=R_2=0$.

Once we have a mesh, we construct a set of basis functions based on the mesh,
such as the piecewise linear functions  $(i=1,2,\cdots, n-1)$

\ \

\begin{minipage}[t]{3.0in}
\vspace{-3.0cm}
\eqmno
\phi_i(x) &=& \left \{ \begin{array}{ll}
\dsp \fr{x-x_{i-1}}{h}     &   \mbox{if $x_{i-1} \le x < x_i$, }\\ \eqsp
\dsp \fr{x_{i+1} - x}{h}     &    \mbox{ if $x_i \le x < x_{i+1}$, }  \\ \eqsp
0 &    \mbox{ otherwise}\, ,
  \end{array} \right.
\enmno
\end{minipage} $\quad$
\begin{minipage}[t]{1.0in}
 \includegraphics[width=1.5\textwidth]{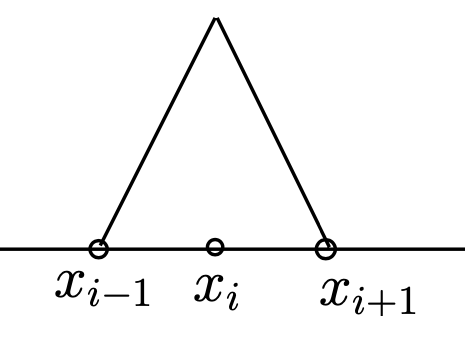}
\end{minipage}

\noindent often called the hat functions, see the right diagram for a hat function.

In the literature, we denote all the piecewise linear functions over the mesh as a finite element space $V_h$,
\eqm
    V_h = \Big \{ v_h(x), \quad  v_h(x) = \sum_{j=1}^{M} c_j \phi_j(x)\, , \quad c_j \in {\cal{R}},  \quad j=1,2,\cdots M \Big\},
\enm
where $M$ is $n-1$, or $n$, or $n+1$, depending on the boundary conditions.
The finite element solution then is a special combination of the basis functions
\eqm
u_h(x) = \sum_{j=1}^{M} c_j^* \phi_j(x)\, ,
\enm
that minimizes the energy error in the space $V_h$.
The coefficients $c_j^*$'s are determined from the  linear system of equations $A {\bf U} = {\bf F} $ with
\eqm
  A = \{ a_{ij}\}, \qquad  a_{ij} = a(\phi_i(x), \phi_j(x)), \qquad   {\bf F} = \{ F_{i}\} , \quad F_i = L(\phi_i(x)),
\enm
 where $ {\bf U} = [ c_1, \, c_2, \,  \cdots,  c_M ]^T$.

 It is well know that the finite element solution converges the true solution as $h\goto 0$ according to the following error estimates
 \eqm
   \| u - u_h\|_{L^2} \le C h^2, \qquad  \| u - u_h\|_{H^1} \le C h, \qquad  \| u - u_h\|_{e } \le C h,
 \enm
 where $C$ is a generic error constant. Thus the finite element method is second order accurate in the $L^2$ norm and first order accurate in the $H^1$ norm.

\section{A new third order method based on a posterior error analysis}  \label{posterior}

Now we consider the special case that $\beta(x)= \beta$ a constant and $q(x) =0$. In this case, the ODE is simply $u''(x) = - f(x)/\beta$.
It is well-known that the finite element solution is the exact at nodal points for the special case.
Thus the interpolation function is the same as the finite element solution,
\eqm
  u_h^ I(x) = \sum_{j=1}^{M} u(x_j)  \phi_j(x)\ = \sum_{j=1}^{M}   c_j^*  \phi_j(x) = u_h(x) ,
\enm
%
From the classical interpolation theory, see for example, \cite{ }, we know that on an element $e_k=[x_k,x_{k+1}]$, the following is true
\eqml{liner-interp}
u(x) &=& \dsp u_h^{I}(x) + \frac{1}{2} \big(x-x_k \big) \big(x-x_{k+1} \big)u''(x) + O(h^3) \\ \eqsp
   &=& \dsp  u_h(x) - \frac{1}{2} \big(x-x_k \big) \big(x-x_{k+1} \big) \, \frac{f(x)}{\beta}  + O(h^3) .
\enml
Thus we obtain a third method using a posterior error estimate with a simple correction term that can be easily calculated.

\begin{theorem}
 Let $u(x)\in H^2(x_l,x_r)$, $f(x) \in L^2(x_l,x_r)$, and $u_h(x)$ be the finite element solution obtained using the $P_1$ finite element space, then
 \eqm
   u_h^{new} (x) = u_h(x) - \frac{1}{2\beta } \big(x-x_k \big) \big(x-x_{k+1} \big) \, f(x), \qquad x_k < x < x_{k+1},
 \enm
 is a third order approximation to the true solution $u(x)$ and we have the following error estimates
 \eqm
   \quad \| u - u_h^{new}\|_{L^2} \le C h^3, \qquad  \| u - u_h^{new}\|_{H^1} \le C h^2, \qquad  \| u - u_h^{new}\|_{e } \le C h^2.
 \enm
 Note also that $ u_h^{new} (x_k)  = u(x_k)$ for all $k$'s.
 That is, the solution has one  more  order accuracy.
\end{theorem}

\section{A new third order compact FE method for variable coefficients}

The above third order method based on a posterior error analysis is valid only if $u_h(x_k) = u(x_k)$ for all $k$'s, which is not be true  in general if $\beta(x)$ is  not a constant. Note that, the third second order method by adding auxiliary points is not compact. In this section, we propose a new compact finite element method based on the interpolation function to construct a bubble function with compact support.

The main idea is to use the interpolation theory and
\eqm
  u''(x) =  - \frac{\beta'(x)}{\beta(x)} u' (x)+ \frac{q(x)}{\beta(x)} u(x) - \frac{f(x)}{\beta(x)} ,
\enm
to get a third order approximation to $u(x)$ and then use  new estimate to construct the basis functions and the FEM solution.
We know that
\eqmno
  u(x) &= & u_h^{I}(x) + \frac{1}{2} u''(x)\big(x-x_k \big) \big(x-x_{k+1} \big) + O(h^3) \\ \eqsp
    &= & u_h^{I}(x) - \frac{1}{2} \frac{\beta'(x)}{\beta(x)} \big(x-x_k \big) \big(x-x_{k+1} \big) u'_{I}(x)  +
 \frac{1}{2}\frac{q(x)}{\beta(x)} \big(x-x_k \big) \big(x-x_{k+1} \big) u_{I}(x) \\ \eqsp
 && \null \qquad-  \frac{1}{2}\frac{f(x)}{\beta(x)} \big(x-x_k \big) \big(x-x_{k+1} \big) + O(h^3) .
\enmno
Note that the added terms has compact support.

The new modified  finite element solution is defined as
\eqml{femnew}
 {u}_h^{new} (x) &=& \dsp  u_{h}(x) - \frac{1}{2} \frac{\beta'(x)}{\beta(x)} \big(x-x_k \big) \big(x-x_{k+1} \big) u'_{h}(x)   \\ \eqsp
&& \dsp  \null +   \frac{1}{2}\frac{q(x)}{\beta(x)} \big(x-x_k \big) \big(x-x_{k+1} \big) u_{h}(x) -  \frac{1}{2}\frac{f(x)}{\beta(x)} \big(x-x_k \big) \big(x-x_{k+1} \big),
\enml
where $u_h(x)$ is defined as before $u_h(x) = \sum_{j=1}^{M} c_j \phi_j(x)$. We use the weak form again to determine the linear  coefficient $c_j^*$.

\begin{theorem}
 Let  $u(x)\in H^3(x_l,x_r)$, $f(x) \in H^1(x_l,x_r)$, $\beta(x)\in C^1(x_l,x_r)\cap H^2(x_l,x_r)$, $ q(x)\in C(x_l,x_r) \cap H^1(x_l,x_r)$, and ${u}_h^{new} (x)$ be the finite element solution obtained using the above formulation, then
 ${u}_h^{new} (x)$  is a third order approximation to the true solution $u(x)$ with the following   error estimates
 \eqm
   \quad \| u - u_h^{new}\|_{L^2} \le C h^3, \qquad  \| u - u_h^{new}\|_{H^1} \le C h^2, \qquad  \| u - u_h^{new}\|_{e } \le C h^2.
 \enm
\end{theorem}

{\bf Proof:} First, we define the new basis function in the interval $(x_k, x_{k+1})$ as
\eqml{femnew2}
  \phi^{new}_i(x)  &=& \dsp  \phi_i(x) - \frac{1}{2} \frac{\beta'(x)}{\beta(x)} \big(x-x_k \big) \big(x-x_{k+1} \big) \phi'_i(x)    \\ \eqsp
&& \dsp  \null +   \frac{1}{2}\frac{q(x)}{\beta(x)} \big(x-x_k \big) \big(x-x_{k+1} \big) \phi_i(x) -  \frac{1}{2}\frac{f(x)}{\beta(x)} \big(x-x_k \big) \big(x-x_{k+1} \big),
\enml
where $\phi_i(x) $ is  an  original hat function. Note that the new basis function depends on the $\beta(x), \, q(x)$ and $f(x)$.
We enlarge the piecewise linear space over the mesh to the new   space below  
\eqm
    V^{new}_h = \Big \{ v^{new}_h(x), \quad  v^{new}_h(x) = \sum_{j=1}^{M} c_j \phi^{new}_j(x)\, , \quad c_j \in {\cal{R}},  \quad j=1,2,\cdots M \Big\}.
\enm
Thus the degree of the freedom of $V^{new}_h $ is the same as that of $V_h$. If we take a special element in $V^{new}_h $ in $(x_k, x_{k+1})$ as
\eqmno
  v_h^{new}(x) &=& \dsp  u^I_{h}(x) - \frac{1}{2} \frac{\beta'(x)}{\beta(x)} \big(x-x_k \big) \big(x-x_{k+1} \big) \frac{d }{d x} u^I_{h}(x)   \\ \eqsp
&& \dsp  \null +   \frac{1}{2}\frac{q(x)}{\beta(x)} \big(x-x_k \big) \big(x-x_{k+1} \big) u^I_{h}(x) -  \frac{1}{2}\frac{f(x)}{\beta(x)} \big(x-x_k \big) \big(x-x_{k+1} \big),
\enmno
where $u^I_{h}(x)$ is the linear interpolation function. Thus we have
\eqmno
  \quad \| u - v_h^{new}\|_{L^2} \le C h^3, \qquad  \| u - v_h^{new}\|_{H^1} \le C h^2, \qquad  \| u - v_h^{new}\|_{e } \le C h^2.
 \enmno
Since the finite element solution is the energy norm and the equivalence of the energy norm and the $H^1$ norm, the theorem  follows directly.

\begin{remark}
  The new method and the convergence theorem are also valid when $\beta(x)$ is a constant. So the discussion in Section \ref{posterior} is included.
  Nevertheless , the algorithms are different, one is a posterior error estimate while the other is a new finite element method.
\end{remark}

\section{Numerical  experiments} \label{sec:result}

In this section, we show two numerical experiments for the one-dimensional Sturm-Liouville problem. We present the $L^2$ and $H^1$ errors.
The order of convergence is estimated using the following formula
 \[
Order= \bigg| \frac{\log(\|E_{N_1}\|_{\infty} /\|E_{N_2} \|_{\infty} }{\log(N_2/N_1)})\bigg|.
\]
with two different $N$'s.

 \begin{example} \label{ex1}
\end{example}
In this example, we show a example for 1D Poisson problem, with the following analytic solution
\begin{equation}
u(x)= \sin(k x).
\end{equation}

In Table~\ref{tab1}-\ref{tab2}, we show grid refinement results for the problem with $k=5\pi$ and $k=50 \pi$, respectively.  In these tables, the first column is the mesh size, the second column is the infinity error of the computed solution and the computed convergence order on the right, and the third is the error and convergence order of the derivative of the solution, respectively.  In all of the cases, we see a clean third order convergence in the solution and second order convergence in the derivative.
\begin{table}[htbp]
\caption{A grid refinement analysis of the third order compact finite element method for Example~\ref{ex1} with $k=5\pi$.} \label{tab1}
\begin{center}
\begin{tabular}{|c|c c|c c|}
\hline
$N$ &         $\|E\|_{L^2}$  & Order  &     $\|E\|_{H^1}$   & Order  \\
\hline
8 	&	9.4674E-01	&		&	1.5922E+01	&		\\
16 	&	7.6248E-02	&	3.63 	&	2.8671E+00	&	2.47 	\\
32 	&	1.1364E-02	&	2.75 	&	8.0230E-01	&	1.84 	\\
64 	&	1.4521E-03	&	2.97 	&	2.0772E-01	&	1.95 	\\
128 	&	1.8229E-04	&	2.99 	&	5.2407E-02	&	1.99 	\\
256 	&	2.2810E-05	&	3.00 	&	1.3132E-02	&	2.00 	\\
512 	&	2.8532E-06	&	3.00 	&	3.2849E-03	&	2.00 	\\
1024 &	3.7524E-07	&	2.93 	&	8.2137E-04	&	2.00 	\\
\hline
\end{tabular}
\end{center}
\end{table}

\begin{table}[htbp]
\caption{A grid refinement analysis of the third order compact finite element method for Example~\ref{ex1} with $k=50\pi$.} \label{tab2}
\begin{center}
\begin{tabular}{|c|c c|c c|}
\hline
$N$ &         $\|E\|_{L^2}$  & Order  &    $\|E\|_{H^1}$  & Order  \\
\hline
64 	&	7.7242E+00	&		&	3.8035E+02	&		\\
128 	&	1.2320E-01	&	5.97 	&	4.4494E+01	&	3.10 	\\
256 	&	2.1751E-02	&	2.50 	&	1.2243E+01	&	1.86 	\\
512 	&	2.8265E-03	&	2.94 	&	3.2238E+00	&	1.93 	\\
1024 &	3.5579E-04	&	2.99 	&	8.1745E-01	&	1.98 	\\
2048 &	4.4562E-05	&	3.00 	&	2.0510E-01	&	1.99 	\\
4096 &	5.8690E-06	&	2.92 	&	5.1322E-02	&	2.00 	\\
\hline
\end{tabular}
\end{center}
\end{table}

 \begin{example} \label{ex2}
\end{example}
In this example, we show a example for 1D elliptic problem, with the following analytic solution
\begin{equation}
u(x)= \sin(k_1 x) \cos(k_2 x).
\end{equation}
The coefficient $\beta (x)$ and $q(x)$ are chosen as
\begin{equation}
\beta(x)= e^{x}, \ \ q(x) =  x^2.
\end{equation}

In Table~\ref{tab3}-\ref{tab6}, we show grid refinement results for the problem with different $k_1$ and $k_2$.  In all of the cases, we see a clean third order convergence in the solution and second order convergence in the derivative.

\begin{table}[htbp]
\caption{A grid refinement analysis of the third order compact finite element method for Example~\ref{ex2} with $k_1 = 5 \pi$ and $k_2 = 0$.} \label{tab3}
\begin{center}
\begin{tabular}{|c|c c|c c|}
\hline
$N$ &         $\|E\|_{L^2}$  & Order  &    $\|E\|_{H^1}$  & Order  \\
\hline
8	&	1.5393E+00	&		&	1.7273E+01	&		\\
16	&	9.5685E-02	&	4.01 	&	3.0023E+00	&	2.52 	\\
32	&	1.2644E-02	&	2.92 	&	8.0654E-01	&	1.90 	\\
64	&	1.5437E-03	&	3.03 	&	2.0909E-01	&	1.95 	\\
128	&	1.8828E-04	&	3.04 	&	5.2686E-02	&	1.99 	\\
256	&	2.3235E-05	&	3.02 	&	1.3196E-02	&	2.00 	\\
512	&	2.8846E-06	&	3.01 	&	3.3004E-03	&	2.00 	\\
1024	&	3.6254E-07	&	2.99 	&	8.2518E-04	&	2.00 	\\
\hline
\end{tabular}
\end{center}
\end{table}

\begin{table}[htbp]
\caption{A grid refinement analysis of the third order compact finite element method for Example~\ref{ex2} with $k_1 = 50 \pi$ and $k_2 = 0$.} \label{tab4}
\begin{center}
\begin{tabular}{|c|c c|c c|}
\hline
$N$ &         $\|E\|_{L^2}$  & Order  &    $\|E\|_{H^1}$  & Order  \\
\hline
64	&	7.9628E+00	&		&	3.8417E+02	&		\\
128	&	1.2781E-01	&	5.96 	&	4.4780E+01	&	3.10 	\\
256	&	2.2066E-02	&	2.53 	&	1.2249E+01	&	1.87 	\\
512	&	2.8471E-03	&	2.95 	&	3.2242E+00	&	1.93 	\\
1024	&	3.5711E-04	&	3.00 	&	8.1755E-01	&	1.98 	\\
2048	&	4.4621E-05	&	3.00 	&	2.0512E-01	&	1.99 	\\
4096	&	5.6445E-06	&	2.98 	&	5.1325E-02	&	2.00 	\\
\hline
\end{tabular}
\end{center}
\end{table}

\begin{table}[htbp]
\caption{A grid refinement analysis of the third order compact finite element method for Example~\ref{ex2} with $k_1 = 5 \pi$ and $k_2 = 5 \pi$.} \label{tab5}
\begin{center}
\begin{tabular}{|c|c c|c c|}
\hline
$N$ &         $\|E\|_{L^2}$  & Order  &    $\|E\|_{H^1}$  & Order  \\
\hline
8	&	9.5172E+00	&		&	1.2581E+02	&		\\
16	&	6.2881E-01	&	3.92 	&	1.6649E+01	&	2.92 	\\
32	&	4.3173E-02	&	3.86 	&	2.9322E+00	&	2.51 	\\
64	&	6.0059E-03	&	2.85 	&	8.0410E-01	&	1.87 	\\
128	&	7.4785E-04	&	3.01 	&	2.0788E-01	&	1.95 	\\
256	&	9.2580E-05	&	3.01 	&	5.2486E-02	&	1.99 	\\
512	&	1.1502E-05	&	3.01 	&	1.3149E-02	&	2.00 	\\
1024	&	1.4323E-06	&	3.01 	&	3.2889E-03	&	2.00 	\\
\hline
\end{tabular}
\end{center}
\end{table}

\begin{table}[htbp]
\caption{A grid refinement analysis of the third order compact finite element method for Example~\ref{ex2} with $k_1 = 50 \pi$ and $k_2 = 50 \pi$.} \label{tab6}
\begin{center}
\begin{tabular}{|c|c c|c c|}
\hline
$N$ &         $\|E\|_{L^2}$  & Order  &    $\|E\|_{H^1}$  & Order  \\
\hline
128	&	3.9235E+00	&		&	3.8229E+02	&		\\
256	&	6.2757E-02	&	5.97 	&	4.4637E+01	&	3.10 	\\
512	&	1.0954E-02	&	2.52 	&	1.2245E+01	&	1.87 	\\
1024	&	1.4185E-03	&	2.95 	&	3.2238E+00	&	1.93 	\\
2048	&	1.7822E-04	&	2.99 	&	8.1748E-01	&	1.98 	\\
4096	&	2.2338E-05	&	3.00 	&	2.0511E-01	&	1.99 	\\
8192&	3.0299E-06	&	2.90 	&	5.1324E-02	&	2.06 	\\
\hline
\end{tabular}
\end{center}
\end{table}

\section{Conclusions and acknowledgements}

A new third order compact finite element method is proposed in this paper. For constant coefficient, the method is a simple posterior error estimate technique and the third order convergence has been proved.
For variable coefficients, the new third order compact finite element method has been proposed and confirmed numerically.   Rigorous proof of the optimal convergence  is  also presented. The degree of freedom and the structure of the resulting linear system of equations of the new third order compact finite element method  are  the same as the traditional finite element method  using piecewise linear functions over the mesh.

Z. Li is partially supported by a Simon's grant 633724.  B. Dong are partially supported by the CNSF Grant No. 11961054.


\parskip= 0.0cm
\renewcommand{\baselinestretch}{1.0}

\bibliographystyle{amsplain}

\bibliography{../../BIB/bib,../../BIB/zhilin,../../BIB/other,anis,referee}






\end{document}